\input amstex \documentstyle{amsppt}
\magnification=1100 \overfullrule=0pt
{\catcode`\@=11\gdef\logo@{}} \hoffset = -0.7 cm
\hsize=6.2truein \vsize=8.5truein \parskip=6pt
plus 2 pt \TagsOnRight

\baselineskip=14pt

\hcorrection{0.5 in}

\def\w{\omega}
\def\P{Proposition 3.}
\def\C{Corollary 3.}
\def\today{\ifcase\month\or
  January\or February\or March\or April\or May\or June\or
  July\or August\or September\or October\or November\or December\fi
  \space\number\day, \number\year
  }

\topmatter
\title MINIMUM $K_{2,3}$-Saturated GRAPHS
\endtitle
\author Ya-Chen Chen \endauthor
\author Ya-Chen Chen \endauthor
\email Ya-Chen.Chen{\@}asu.edu
\endemail

\abstract  A graph is {\it $K_{2,3}$-saturated} if it has
no subgraph isomorphic to $K_{2,3}$, but does contain a 
$K_{2,3}$ after the addition of any new edge.
We prove that the minimum number of edges in a
$K_{2,3}$-saturated graph on $n \geq 5$ vertices is 
$sat(n, K_{2,3}) = 2n - 3$.
\endabstract
${}$
\date  \today
\footnote"*"{
 This copy was printed on \today}
\enddate
\endtopmatter
\document

\subhead \S\ 1. Introduction \endsubhead
\medskip

All graphs studied are simple ones. 
We denote a path, a cycle, a star, a complete 
graph, the complement of a complete graph, and a 
complete $r$-uniform hypergraph with $n$ vertices 
by $P_n$, $C_n$, $S_n$, $K_n$, $I_n$, and $K_n^r$, respectively. 
We write $K_{n_1,...,n_r}$ for
the complete $r$-partite graph with partite sets
of sizes $n_1,...,n_r$.
For a graph $G$, denote $V = V(G)$, let $N(x)$ be the set of the
vertices adjacent to $x$, and $d(x) = |N(x)|$, $N[x] = N(x) \cup
\{ x \}$, $e(G) = |E(G)|$, and $\delta(G) = \min \{ d(x): x \in 
V \}$.
If $A, B \subseteq V$, we define $G[A,B]$ to be the
subgraph with vertex set $A \cup B$ and edge set
$E(G[A,B]) = \{ xy \in E(G): x \in A, y \in B \}$.
We write $e(G[A,B])$ as $e([A,B])$. If $A = B$, we write $G[A,
A]$ as $G[A]$ and $|E(G[A])|$ as $e(A)$.
A {\it clique} in a
graph is a set of pairwise adjacent vertices.
For two graphs $G$ and $H$, the {\it disjoint union} $G + H$
has vertex set $V(G) \cup V(H)$ and edge set $E(G) \cup E(H)$.
The {\it join} $G*H$ is obtained from $G + H$ by adding the 
edges $\{ xy: x \in V(G), y \in V(H) \}$.

Let $\Cal F$ be a family of graphs or hypergraphs.
A hypergraph is {\it ${\Cal F}$-saturated} if it 
has no $F \in {\Cal F}$ as a subhypergraph, but 
does contain some $H \in {\Cal F}$ after the 
addition of any new edge. 
The minimum and maximum number of edges 
in an ${\Cal F}$-saturated graph is
denoted by $sat(n, {\Cal F})$ and $ex(n, {\Cal F})$,
respectively. 
An ${\Cal F}$-saturated graph $G$ on $n$ vertices 
with $e(G) = sat(n, {\Cal F})$ is called a 
$sat(n, {\Cal F})$-graph. The problem
of determining $ex(n, {\Cal F})$ is Tur\'an's problem. 
If ${\Cal F} = \{ F \}$, we 
also write $sat(n, {\Cal F})$ as $sat(n, F)$.
Erd\H{o}s, Hajnal, and Moon~[9] proved that the 
$sat(n, K_k)$-graph is $K_{k-2} * I_{n-k+2}$.
K\'aszonyi and Tuza [15] determined 
$sat(n, F)$ for $F = S_k, kK_2$, $P_k$, and they 
proved that $sat(n, {\Cal F}) = O(n)$ for any family 
$\Cal F$ of graphs. 

As for hypergraphs, Bollob\'as [4] 
generalized Erd\H{o}s, Hajnal, and Moon's results 
to $K_k^r$-saturated hypergraphs.
Erd\H{o}s, F\"uredi, and Tuza [8] obtained $sat(n, 
F)$ for some particular hypergraphs $F$ 
with few edges. 
Pikhurko [17] proved Tuza's 
conjecture that $sat(n, {\Cal F}) = O(n^{r-1})$ 
for all families of $r$-uniform hypergraphs 
whose independence numbers are bounded by a 
constant. For more results and open problems, see [18].

For cycles, Ollmann [16] pointed out that the $sat(n,       
C_3)$-graph is the star $S_n$, and he obtained all   
$sat(n, C_4)$-graphs. Later Tuza [21] gave a shorter
proof for $sat(n, C_4)$. Ashkenazi [1] described the
properties of $C_3$-saturated graphs, planar
$C_3$-saturated graphs, and $C_4$-free $C_3$-saturated
graphs. Fisher et al.
[11] constructed the $C_5$-saturated graphs establishing
the upper bounds of $sat(n, C_5)$. Chen [6,7] determined all 
$sat(n, C_5)$-graphs. Barefoot et al.
[2] showed that $n + c_1n/k \leq sat(n, 
C_k) \leq n + c_2n/k$ for some positive $c_1, c_2$.
They [2] and Gould, Luczak, and Schmitt [13]
gave new upper bounds for $sat(n, C_k)$ for small $k$.
Recently, F\"uredi and Kim [12] gave almost exact asymptotics 
for $sat(n,C_k)$ as $k$ is fixed and $n \rightarrow \infty$.
For more results and open problems, see the excellent
survey by Faudree, Faudree, and Schmitt [10].

Pikhurko [17] and G. Chen et al. [5]
obtained $sat(n, K_{1,...,1,\ell})$ of the 
complete $(r+1)$-partite graph $K_{1,...,1,\ell}$ for $n \geq 
n(r, \ell)$. 
Pikhurko and 
Schmitt [19] presented $K_{2,3}$-saturated graphs with 
$2n - 3$ edges and proved $sat(n, K_{2,3}) \geq 2n - cn^{3/4}$, 
where $c$ is a constant. Gould and Schmitt [14] conjectured
that the complete $r$-partite graph $K_{2,...,2}$ has
$sat(n, K_{2,...,2}) = \lceil ((4r-5)n -4r^2 + 6r - 1)/2 
\rceil$ and proved it when the minimum degree of the 
$K_{2,...,2}$-saturated graphs is $2r - 3$. 
Recently, Bohman, Fonoberova, and Pikhurko [3] proved 
that for $r \geq 2$ and $s_r \geq ... \geq s_1 \geq 1$,
as $n \rightarrow \infty$, 
$sat(n, K_{s_1,...,s_r}) = (s_1 + ... + s_{r-1} + 0.5s_r - 
1.5)n + O(n^{3/4}).$ 
They [3] constructed a
$K_{s_1,...,s_r}$-saturated graph $K_p * H$ with $(s_1 + ... 
+ s_{r-1} + 0.5s_r - 1.5)n + O(1)$ edges, where $H$ is
a $K_{1,s_r}$-saturated graph and
$p = s_1 + ... + s_{r-1} - 1$. They [3] showed
that any $K_{s_1,...,s_r}$-saturated graph on $n$ vertices
with at most $sat(n, K_{s_1,...,s_r}) + o(n)$ edges can be 
transformed into $K_p * H$ by adding and removing at most
$o(n)$ edges. Bohman, Fonoberova, and Pikhurko [3] also
conjectured that $sat(n, K_{2,3}) = 2n - 3$. Here we prove 
their conjecture.

\proclaim {Theorem 1} 
$sat(n, K_{2,3}) = 2n - 3$.
\endproclaim

We present $sat(n, K_{2,3})$-graphs in Section 2,
obtain the properties and structures of 
$K_{2,3}$-saturated graphs in Section 3, and prove
Theorem 1 in Sections 4-6.
\subhead \S\ 2. Extremal graphs \endsubhead
\medskip

In this section, we present $sat(n, K_{2,3})$-graphs.  
The construction of these graphs was applied to obtain
general upper bounds of $sat(n, {\Cal F})$ for any family
${\Cal F}$ of graphs by K\'aszonyi and Tuza [15], and 
the structure was recently
proved to be possessed by all almost extremal 
$K_{s_1,...,s_r}$-saturated graphs by Bohman, Fonoberova, and 
Pikhurko [3]. Pikhurko and Schmitt [19] have 
presented these $K_{2,3}$-saturated graphs with $2n - 3$ edges;
we include them for completeness. 
Let ${\Cal H}_n$ be the
set of $2$-regular $K_{2,2}$-free graphs on $n$ vertices. 
Let ${\Cal G}_n =$ 
$\{ K_1 * (K_1 + H): H \in {\Cal H}_{n-2}
\}$ $\cup \{ K_1 * (K_2 + H): H \in {\Cal H}_{n-3} \} \cup
\{ K_1 * (P_4 +  H): H \in {\Cal H}_{n-5} \}.$ 
It is easy to verify that each of the 
graphs in ${\Cal G}_n$ is $K_{2,3}$-saturated, since the 
addition of any edge results in either a $K_{1,3}$ or a
$C_4$ disjoint from the vertex of maximum degree.

\subhead \S\ 3. Properties of $K_{2,3}$-saturated graphs
\endsubhead
\medskip 

In this section, let $G$ be a $K_{2,3}$-saturated 
graph and we shall describe properties of $G$. 

\proclaim{\P1} If $\alpha_1 \alpha_2 \not \in E(G)$, then there
exists a vertex $b$ in $N(\alpha_i)$ such that $b$ and 
$\alpha_{3-i}$ have two common neighbors, namely, 
$|N(\alpha_{3-i}) \cap N(b)| \geq 2$, $i = 1$ or $2$.
\endproclaim
\demo{Proof} Let $A, B$ be two partite sets generating a
$K_{2,3}$ obtained by adding $\alpha_1 \alpha_2$ to $E(G)$.
Let $A = \{ \alpha_i, z, w \}$ and $B = \{ 
\alpha_{3-i}, b \}$.
Thus $b \in N(\alpha_i)$ is adjacent to
$z, w \in N(\alpha_{3-i})$.
\qed\enddemo

\proclaim{\C2} Let $N(\alpha) = \{ x_1, .., x_k \}$. If
$y \in V \setminus N[\alpha]$, then either there exists $l$ such 
that $|N(y) \cap N(x_l)| \geq 2$ or there are $i, j$
such that $i \neq j$ and $y$ is adjacent to 
a vertex in $N(x_i) \cap N(x_j)$. 
In particular, if $k = 1$, then for all $y \in V \setminus 
N[\alpha]$, it satisfies $|N(y) \cap N(x_1)| \geq 2$.
\endproclaim 

\proclaim{\C3} Let $V_1 = N[\alpha] \cup \{ v \in V: |N(v) 
\cap N(\alpha)| = 2 \}$, the set ${\Cal U}_2 = \{ v \in V 
\setminus V_1:$ $|N(v) \cap N(\alpha)| = 1 \},$ and ${\Cal 
U}_3 = V \setminus (V_1 \cup {\Cal U}_2)$.
We define $\w(b) = |N(b) \cap V_1| + 0.5|N(b) \cap {\Cal 
U}_2|$ for $b \in {\Cal U}_2$.
Let $x^* \in N(\alpha)$ with $|N(x^*) \cap N(\alpha)| \leq 1$.
When $y \in {\Cal U}_2$ and  
$x^*$ is the unique common neighbor of $\alpha$ and $y$,
we conclude $\w(y) \geq 1.5$ and if $\w(y) = 1.5$, then 
there exist $x \in N(\alpha) \cap N(x^*)$ and $y' \in N(x)$ 
such that $N(y) \cap {\Cal U}_2 = \{ y' \}$. 
If $z \in {\Cal U}_3$, then we have $|N(z) \cap (V \setminus 
{\Cal U}_3)| \geq 1 + |\{ x \in N(\alpha): N(z) \cap N(x) 
\cap {\Cal U}_2 \neq \emptyset \}|.$
\endproclaim 
\demo{Proof} 
We assume $N(y) \cap V_1 = \{ x^* \}$.
Since $|N(x^*) \cap N(\alpha)| \leq 1$, by \C2, there is $x \in
N(\alpha)$ such that $|N(y) \cap N(x)| \geq 2$. Thus there is $y'
\in (N(y) \cap N(x)) \setminus \{ x^* \}$ and $\w(y) \geq 1.5$.
If $\w(y) = 1.5$, then $N(y) \cap {\Cal U}_2 = \{ y' \}$ and $N(y)
\cap N(x) = \{ y', x^* \}$.
Since $G$ is $K_{2,3}$-free, for $z \in V(G) \setminus \{ \alpha
\}$, $|N(z) \cap N(\alpha)| \leq 2$. Let $z \in {\Cal U}_3$. Thus
$z \alpha \not \in E(G)$ and $|N(z) \cap N(\alpha)| = 0$. 
Since every
vertex in $N(z) \cap {\Cal U}_2$ is adjcent to exactly one vertex
in $N(\alpha)$, it follows $|N(z) \cap {\Cal U}_2| \geq |\{ x \in 
N(\alpha): N(z) \cap N(x) \cap {\Cal U}_2 \neq \emptyset \}|$. 
Let $|N(z) \cap V_1| = 0$. Thus every
vertex in $N(z)$ is adjacent to at most one veretx in $N(\alpha)$.
By \P1, there is $x_1 \in N(\alpha)$ adjacent 
to two vertices $w_1, w_2 \in N(z)$. Since every vertex 
in $(N(z) \setminus \{ w_1 \}) \cap {\Cal U}_2$ is adjcent to 
exactly one vertex in $N(\alpha)$, the vertex $w_1$ is adjacent 
to only $x_1$ in $N(\alpha)$, and $x_1$ has $w_2 \in N(x_1) 
\cap N(z) \cap {\Cal U}_2$, it follows $|N(z) \cap (V \setminus 
{\Cal U}_3)| = |\{ w_1 \}| + |(N(z) \setminus \{ w_1 \}) \cap 
{\Cal U}_2| \geq 1 + |\{ x \in N(\alpha): N(z) \cap N(x) \cap 
{\Cal U}_2 \neq \emptyset \}|$.
\qed\enddemo

\bigskip 
\subhead \S\ 4. The case when $\delta(G) = 1$
\endsubhead
\medskip
In this section, we prove Theorem 1 when $\delta(G) = 1$.
Although Pikhurko and Schmitt [19] have proved Theorem 1 
when $\delta(G) = 1$, for completeness,
we include our proof here.

\proclaim{Lemma 4.1} If $G$ is a sat$(n, K_{2,3})$-graph
with $\delta(G) = 1$, then $e(G) \geq 2n - 3$.
\endproclaim 
\demo{Proof} Let $N(\alpha) = \{ x \}$, $U_1 = N(x) \setminus 
\{ \alpha \}$, and $U_2 = V \setminus N[x]$. By \C2, if $y \in 
U_1 \cup U_2$, then $|N(y) \cap U_1| \geq 2$. Hence
$$e(G) = d(x) + e(U_1) + e([U_1, U_2]) + e(U_2)
= d(x) + 0.5\sum_{y \in U_1} |N(y) \cap U_1| + 
\sum_{z \in U_2} |N(z) \cap U_1| + e(U_2)$$
$\geq d(x) + |U_1| + 2|U_2| + e(U_2) = d(x) + (d(x) - 1) + 2(n - d(x) - 1)
+ e(U_2) \geq 2n - 3 + e(U_2)$. 
\qed
\enddemo

\bigskip
\subhead \S\ 5. The case when $\delta(G) = 2$
\endsubhead
\medskip
In this section, we prove Theorem 1 when $\delta(G) = 2$.

\proclaim{Lemma 5.1} Let $G$ be a sat$(n, K_{2,3})$-graph with 
$\delta(G) = 2$, $A$ be the set of degree $2$ vertices having adjacent
neighbors, and $B$ be the set of degree $2$ vertices whose neighbors have
exactly one common neighbor.
If $A \cap B \neq \emptyset$, then $e(G) \geq 2n 
- 3$.
\endproclaim
\demo{Proof} We choose a vertex $\alpha \in A \cap B$ and denote 
$N(\alpha) = \{ x_1, x_2 \}$. We define $V_1 =  N[\alpha]$, ${\Cal U}_2 
= (N(x_1) \cup N(x_2)) \setminus V_1$, and ${\Cal U}_3 = V \setminus
(V_1 \cup {\Cal U}_2)$. For $y \in {\Cal U}_2$, we define
$$\w(y) = |N(y) \cap V_1| + 0.5|N(y) \cap {\Cal U}_2| - 2,$$
$$U_2^+ = \{ y \in {\Cal U}_2:                 
\w(y) \geq 0.5 \}, ~{\Cal U}_2^- = \{ y \in {\Cal             
U}_2: \w(y) < 0 \}, ~\text{and} ~{\Cal U}_2^0 = ~{\Cal U}_2     
\setminus ~{\Cal U}_2^-.$$
For $z \in {\Cal U}_3$, let $\w(z) = |N(z) \cap
{\Cal U}_2|$, 
$$U_3^3 = \{ z \in {\Cal U}_3: \w(z) = 3 \},
~\text{and} ~U_3^4 = \{ z \in {\Cal U}_3: \w(z) \geq 4 \}.$$
For $y \in {\Cal U}_2^-$, we define $f(y)$, a subset of $N(y)$.
We partition ${\Cal U}_2^-$ into $S_0, ..., S_4$:
$$S_0 = \{ y \in {\Cal U}_2^-: f(y) \subseteq U_2^+ 
\}, ~S_4 = \{ y \in {\Cal U}_2^-: f(y) \cap U_3^4 \neq \emptyset 
\},$$
$$S_1 = \{ y \in {\Cal U}_2^-: |f(y)| = 1, f(y) \subseteq U_3^3 
\}, ~S_2 = \{ y \in {\Cal U}_2^-: |f(y)| = 2, f(y) \subseteq U_3^3
\}, ~\text{and}$$ 
$$S_3 = \{ y \in {\Cal U}_2^-: f(y)
= \{ y', z \}, z \in U_3^3, y' \in {\Cal U}_2 \}.$$
Now we define $f(y)$ and verify that ${\Cal U}_2^- = S_0 \cup 
... \cup S_4$.
Let $y \in {\Cal U}_2^- \cap N(x_i)$. By \C3, $\w(y) = - 0.5$
and there is
$y' \in N(x_{3-i}) \setminus V_1$ 
such that $N(y) \cap (V_1 \cup {\Cal U}_2) = \{ x_i, y' \}$. 
If $y' \in U_2^+$, then let $f(y) =
\{ y' \}$. Since $N(y') \cap N(x_i) = \{ x_{3-i},
y \}$, it follows $f(b) \neq \{ y' \}$ if $b \neq y$. Thus
$$|S_0| \leq |U_2^+|. \tag1$$
Now we define $f(y)$ for $y \in
{\Cal U}_2^-$ whose unique neighbor $y'$ in ${\Cal U}_2$ has
$\w(y') \leq 0$.   
If there is $z \in N(y)$ adjacent to distinct $y_1, y_2 \in N(x_{3-i})$, 
then let $f(y) = \{ z \}$. Since $N(x_1) \cap N(x_2) = \{ \alpha \}$,
and $w(y') \leq 0$, 
and $y, y_1, y_2 \in N(z)$, it follows $f(y) \subseteq U_3^3 \cup 
U_3^4$ and $y \in S_4 \cup S_1$. 
Let $|N(z) \cap N(x_{3-i})| \leq 1$ for all $z \in N(y)$.
Since $y x_{3-i} \not \in E(G)$, by \P1, there is $y_3 \in
N(x_{3-i})$ adjacent to distinct $z_1, z_2 \in N(y)$.   
We define $f(y) = \{ z_1, z_2 \}.$ Since $N(y) \cap {\Cal 
U}_2 = \{ y' \}$, it follows $\{ z_1, z_2 \} \not \subset 
{\Cal U}_2$ and $f(y) \cap {\Cal U}_3 \neq \emptyset$. By 
\C3, if $z_j \in f(y) \cap {\Cal U}_3$ and $|f(y)| = 2$, 
then $z_j \in U_3^3 \cup U_3^4$ and $y \in S_2 \cup S_3 \cup
S_4$.
Recall $|N(z) \cap {\Cal U}_2| \geq 4$ for $z \in U_3^4$.
We observe that 
$$(0.5 + 0.5)\sum_{ z \in U_3^4} \w(z) \geq
2|U_3^4| + 0.5|\{ y \in {\Cal U}_2: N(y) \cap U_3^4 \neq
\emptyset \}| \geq 2|U_3^4| + 0.5|S_4|. \tag2$$
We shall partition $U_3^3$ into $U_{31}^3$, $U_{32}^3$,
$U_{33}^3$, and $U_3^3 \setminus (U_{31}^3 \cup               
U_{32}^3 \cup U_{33}^3)$ depending on the manner in which 
$f(y)$ intersect with $U_3^3$. 
Recall that if $z \in U_3^3$, then $\w(z) = |N(z) 
\cap {\Cal U}_2| = 3$.
If $y \in N(x_i)$, $f(y) = \{ z \}$, and $z \in U_3^3$, then
there are $y_1, y_2$ such that $N(z) \cap {\Cal U}_2 \cap N(x_i) 
= \{ y \}$ and $N(z) \cap {\Cal U}_2 \cap N(x_{3-i}) = \{ y_1, 
y_2 \}$.
Thus if $z \in U_3^3$, then $\{ z \}$ is $f(y)$ for at most one
$y \in S_1$.
$$\text{Let} ~U_{31}^3 = \{ z \in U_3^3 \cap f(y): y \in S_1 \}.
~\text{Thus} ~|S_1| = |U_{31}^3|. \tag3$$
Let $y_3 \in N(x_i) \cap S_3$, $y' \in {\Cal U}_2$,
$z' \in U_3^3$, and $f(y_3) = \{ y', z' \}$. 
Thus there is $y_5 \in N(x_{3-i})$ such that $N(y_3) \cap 
N(y_5) = \{ y', z' \}$.       
By \C3, $y' \in N(x_{3-i})$ and $y_5 \in {\Cal U}_2^0$.
Since $z'$ is adjacent to three vertices in ${\Cal U}_2$, 
one of which is in ${\Cal U}_2^- \cap N(x_i)$, another one of
which is in ${\Cal U}_2^0 \cap N(x_{3-i})$, it follows $z'$ 
is in $f(y_3)$ for at most two $y_3 \in S_3$ and $z' \not 
\in U_{31}^3$. We define
$$U_{32}^3 = \{ z \in U_3^3 \cap
f(y) \cap f(b): y, b \in S_3, y \neq b \}
~\text{and} ~U_{33}^3
= \{ z \in (U_3^3 \cap f(y)) \setminus U_{32}^3: y \in S_3 
\}.$$
Thus $$|S_3| = 2|U_{32}^3| + |U_{33}^3| \tag4$$
and $U_{31}^3, U_{32}^3, U_{33}^3$ are pairwise disjoint.
Also if $z \in U_{32}^3$, then $z$ is adjacent to exactly 
three vertices in ${\Cal U}_2$, one of which is in ${\Cal 
U}_2^0$, two of which are in $S_3$, and $N(z) \cap S_2 = 
\emptyset$. If $z' \in U_{33}^3$, then $z'$ is adjacent to 
one vertex in ${\Cal U}_2^0$, one vertex in $S_3$, and 
$|N(z') \cap S_2| \leq 1$.             
Since each $y \in S_2$ is adjacent to two vertices in $U_3^3$,
and each vertex in $U_3^3$ is adjacent to exactly three vertices
in ${\Cal U}_2$,
it follows $2|S_2| \leq 2|U_{31}^3| + |U_{33}^3| +
3|U_3^3 \setminus (U_{31}^3 \cup U_{32}^3 \cup U_{33}^3)|$ and
by (3), (4),
$$2|U_3^3| = |U_{31}^3| + (2|U_{32}^3| + |U_{33}^3|) +
|U_{31}^3| + |U_{33}^3| + 2|U_3^3 \setminus (U_{31}^3 \cup
U_{32}^3 \cup U_{33}^3)| \geq |S_1| + |S_3| + |S_2|. \tag5$$  
If $b \in {\Cal U}_3$, then by \C3, $\w(b) \geq 2$. 
By (1), (2), and (5),
$$e(G) = e(V_1) + e([V_1, {\Cal U}_2]) + e({\Cal U}_2)
+ e([{\Cal U}_2, {\Cal U}_3]) + e({\Cal U}_3)$$
$$\geq 2|V_1| - 3 + 2|{\Cal U}_2| + 0.5(|U_2^+| - 
|S_0| - ... - |S_4|) + 2|{\Cal U}_3| + |U_3^3| + 0.5|S_4| +  e({\Cal   
U}_3) \geq 2n - 3.\qed$$
\enddemo

\proclaim{Lemma 5.2} Let $G$ be a sat$(n, K_{2,3})$-graph with
$\delta(G) = 2$ and $A$ be the set of degree $2$ vertices having adjacent
neighbors. 
If $A \neq \emptyset$, then $e(G) \geq 2n
- 3$.   
\endproclaim   
\demo{Proof} By Lemma 5.1, 
since $G$ is $K_{2,3}$-free, we can assume that for each
$\alpha \in A$, there is a unique vertex $b$ such that $N(\alpha)
\subseteq N(b)$ and $b \neq \alpha$. Let ${\Cal B} = \{ b \in
V: N(b) \supseteq N(a)$ for a vertex $a \in A$ with $a \neq
b \}$. We choose $\alpha \in A$ satisfying the unique vertex 
$\beta$ with $N(\alpha) \subseteq N(\beta)$ and $\beta \neq 
\alpha$ has $d(\beta) = \min \{ d(b): b \in {\Cal B} \}$. Let 
$N(\alpha) = \{ x_1, x_2 \}$. Let $V_1 = {\Cal U}_1 = N[\alpha] \cup 
\{ \beta \}$, ${\Cal U}_2 = (N(x_1) \cup 
N(x_2)) \setminus V_1$, the set ${\Cal U}_3 = V \setminus (V_1 \cup 
{\Cal U}_2 \cup N(\beta))$, and 
${\Cal U}_4 = V \setminus (V_1 \cup {\Cal U}_2 \cup {\Cal U}_3).$
$$~\text{For} ~i \in \{ 2, 3, 4 \} ~\text{and} ~y \in {\Cal U}_i,
~\text{let} ~\w(y) = |N(y) \cap ({\Cal U}_1 \cup ... \cup {\Cal 
U}_{i-1})| + 0.5|N(y) \cap {\Cal U}_i| - 2. \tag6$$
$$\text{Let} ~U_2^+ = \{ y \in {\Cal U}_2: \w(y) \geq 0.5 \}, 
~\text{the set} ~{\Cal U}_2^- = \{ y \in {\Cal
U}_2: \w(y) < 0 \}, ~\text{and} ~{\Cal U}_2^0 = ~{\Cal U}_2
\setminus ~{\Cal U}_2^-.$$ 
$$\text{Let} ~U_3^0 = \{ z \in {\Cal U}_3: \w(z) = 0.5 \},
~U_3^1 = \{ z \in {\Cal U}_3: 1 \leq \w(z) < 2 \},
~\text{and} ~U_3^2 = \{ z \in {\Cal U}_3: \w(z) \geq 2 \}.$$ 
Since $\alpha \beta \not \in E(G)$, by \P1,
without loss of generality, there is $\gamma \in N(x_1) \cap
N(\beta)$. 
For $y \in {\Cal U}_2^-$, we define $f(y)$, a subset of $N(y)$.
For one particular case we also define $f(\gamma)$ and leave
$f(y)$ undefined for exatcly one $y \in {\Cal U}_2^-$.
We partition ${\Cal U}_2^- \cup \{ \gamma \}$ into $S_0, ..., 
S_6$: 
$$\text{let} ~S_0 = \{ y \in {\Cal U}_2^-: f(y) \subseteq U_2^+, 
f(y) \neq \emptyset \}, ~S_i = \{ y \in {\Cal
U}_2^- \cup \{ \gamma \}: |f(y)| = i, f(y) \subseteq U_3^1
\}, i = 1, 2,$$ 
$$S_4 = \{ y \in {\Cal U}_2^- \cup \{ \gamma \}: f(y) \cap U_3^2 
\neq \emptyset \}, ~S_5 = \{ y \in {\Cal U}_2^- \cup \{ \gamma \}:
f(y) \cap {\Cal U}_4 \neq \emptyset, f(y) \cap U_3^2 = \emptyset 
\},$$                 
$$S_3 = \{ y \in {\Cal U}_2^-: f(y) = \emptyset \},
~\text{and} ~S_6 = \{ y \in {\Cal U}_2^-: f(y) = \{ z \}, z \in 
U_3^0 \}.$$
Note that $S_0, ..., S_6$ are pairwise disjoint.
Now we define $f(y)$, verify that each $f(y)$ is as described in 
the definitions of $S_0, ..., S_6$, $|{\Cal U}_2^-| = |S_0| +   
... + |S_6|$, and $f(y)$ satisfies that
$$d(w) \geq 3 ~\text{for} ~w \in {\Cal U}_4 \cap f(y), ~\text{and 
if} ~y \in S_1 ~\text{and} ~z \in f(y), ~\text{then} ~N(z) \cap 
(({\Cal U}_2^0 \setminus \{ \gamma \}) \cup {\Cal U}_3) \neq 
\emptyset. \tag7$$
Let $y \in {\Cal U}_2^- \cap N(x_i)$. By \C3, $\w(y) = - 0.5$ 
and there is
$y' \in N(x_{3-i}) \setminus V_1$ 
such that $N(y) \cap (V_1 \cup {\Cal U}_2) = \{ x_i, y' \}$. 
If $y' \in U_2^+$, then we define $f(y) =
\{ y' \}$ and $y \in S_0$. Since $N(y') \cap N(x_i) = \{ 
x_{3-i}, y \}$, it follows $f(b) \neq \{ y' \}$ if $b \neq y$. 
Since $N(\gamma) \cap N(x_2) = \{ x_1, \beta \}$, it follows $\{ \gamma \} 
\neq f(y)$ for any $y \in S_0$. After 
defining $f(y)$ for all $y \in {\Cal U}_2^-$ having a neighbor in 
$U_2^+$,
we conclude
$$|S_0| \leq |U_2^+ \setminus \{ \gamma \}|. \tag8$$
Let $y \in {\Cal U}_2^- \setminus S_0$. If there is $z \in 
U_3^2 \cap N(y)$ or $z \in {\Cal U}_4 \cap N(y)$ with $d(z) 
\geq 3$, then we choose one such $z$ and define $f(y) = \{ 
z \}$. We have $y \in S_4 \cup S_5$ and $f(y)$ satisfies (7).

Now we define $f(y)$ for $y \in {\Cal U}_2^- 
\setminus (S_0 \cup S_4 \cup S_5)$. 
First let there be $y_1, y_2 \in 
{\Cal U}_2^- \cap N(x_i)$ whose $f(y_1), f(y_2)$ are to be 
defined. By \C3, there is $y_{j+2} \in N(x_{3-i})$ such that 
$N(y_j) \cap (V_1 \cup {\Cal U}_2) = \{ x_i, y_{j+2} \}$, 
$j = 1, 2$. 
Since $y_1 y_2 \not \in E(G)$, by \P1, without loss of 
generality, there is $b \in N(y_1)$ adjacent to  
$z_1, z_2 \in N(y_2)$. Since ${\Cal U}_4 \subset N(\beta)$, 
$y_1, y_2 \not \in S_4 \cup S_5$, and $y_2 \in {\Cal U}_2^-$,
it follows $b \in \{ y_3 \} \cup {\Cal U}_3$ and $z_1, z_2 
\in \{ y_4 \} \cup {\Cal U}_3$. Let $b = y_3$. We define 
$f(y_2) = \{ z_1, z_2 \} \setminus \{ y_4 \}$.
By \C3, since $y_2 \not \in S_4$, it follows
$f(y_2) \subseteq U_3^1$ and $y_2 \in S_1 \cup S_2$. If 
$f(y_2) = \{ z_j \}$, then $y_4, y_1 \in N(b) \cap {\Cal
U}_2$, the vertex $b \in ({\Cal U}_2^0 \setminus \{ \gamma 
\}) \cap N(z_j)$, and $f(y_2)$ satisfies (7). 
Let $b \in {\Cal U}_3$. We define $f(y_1) = \{ b \}$.
By \C3, $b \in U_3^1$ and $y_1 \in S_1$. 
Since $\{ z_1, z_2 \} \cap {\Cal U}_3 \cap N(b) \neq \emptyset$, 
it follows $f(y_1)$ satisfies (7).
We have defined one of $f(y_1), f(y_2)$. 
We repeat this process until there is at most one $y \in 
{\Cal U}_2^- \cap N(x_i)$ whose $f(y)$ is not defined. 

If there is only one $y \in {\Cal U}_2^-$ whose $f(y)$ has not 
been defined, then we define $f(y) = \emptyset$ and $S_3 = \{ y \}$.
Let $u_i^*$ be the unique vertex in ${\Cal U}_2^- \cap N(x_i)$ whose 
$f(u_i^*)$ is not defined, $i = 1, 2$. Recall $\gamma \in N(\beta) 
\cap N(x_1)$. If $\gamma \in U_2^+$, then let 
$f(u_1^*) = f(u_2^*) = \emptyset$. Thus $S_3 = \{ u_1, u_2 \}$. By (8),
$$|U_2^+| - |S_0| - |S_3| \geq - 1. \tag9$$
Let $\w(\gamma) = 0$ and $N(\gamma) \cap (V_1 \cup {\Cal U}_2)
= \{ x_1, \beta \}$. Let there be $w \in 
N(u_1^*) \cap {\Cal U}_4$. Since $f(u_1^*)$ is not defined, 
$N(w) = \{ u_1^*, \beta \}$. Since $N(u_1^*) \cap N(\beta) = \{ x_1, w \}$
and $u_2^* w \not \in E(G)$, by \P1, 
there is $b \in N(w)$ adjacent to $w_1, w_2 \in N(u_2^*)$. If
$b = \beta$, then $w_1$ or $w_2$, say $w_1 \in N(\beta) \setminus 
\{ x_2 \}$. 
Since $u_2^* \in {\Cal U}_2^-$ and $f(u_2^*)$ is not defined, it follows
$w_1 \not \in {\Cal U}_2$, $w_1 \in {\Cal U}_4$, and $N(w_1) = \{ 
\beta, u_2^* \}$, a contradiction to \P1, since $w w_1 \not \in E(G)$. 
Hence $b = u_1^*$. 
Since ${\Cal U}_4 \subset N(\beta)$, by \C3, $w_1, w_2 \in U_3^1$. 
We define $f(u_1^*) = 
f(u_2^*) = \{ w_1, w_2 \}$ and $u_1^*, u_2^* \in S_2$. 

Let $N(u_1^*) \cap {\Cal U}_4 = \emptyset$. We define $f(u_2^*) = 
\emptyset$ and $u_2^* \in S_3$. We consider the $K_{2,3}$ 
created by adding $u_1^* \gamma$.
Let there be $z \in N(u_1^*)$ adjacent to $w_1, w_2 \in N(\gamma)$. 
Let $z \in {\Cal U}_2$. By \C3, $z \in N(x_2)$ and $w_1, w_2 \in U_3^1 
\cup U_3^2 \cup {\Cal U}_4$. We define $f(\gamma) = \{ w_1, w_2 \}$, 
leave $f(u_1^*)$ undefined, and $\gamma \in S_2 \cup S_4 \cup S_5$.
Note that $|{\Cal U}_2^-| = |S_0| + ... + |S_6|$.
Since ${\Cal U}_4 \subset N(\beta)$,
$$\text{if} ~w \in {\Cal U}_4 \cap f(\gamma), ~\text{then} 
~|N(w) \cap (V \setminus {\Cal U}_4)| \geq 3 ~\text{and} ~f(\gamma)
~\text{satisfies} ~(7). \tag10$$
Next let $z \in {\Cal U}_3$. We define $f(u_1^*) = \{ z \}$.
Since $x_1 \in N(\gamma) \cap N(\beta)$, it follows $\{ w_1,
w_2 \} \cap {\Cal U}_3 \cap N(z) \neq \emptyset$. By \C3,
$u_1^* \in S_6 \cup S_1$ and $f(u_1^*)$ satisfies (7). 
Next let $|N(z) \cap N(\gamma)| \leq 1$ for all $z \in N(u_1^*)$. By \P1,
there is $z_1 \in N(\gamma)$ adjacent to $z_2, z_3 \in N(u_1^*)$. 
Since $N(u_1^*) \cap {\Cal U}_4 = \emptyset$ and $N(\gamma) \cap (V_1 
\cup {\Cal U}_2) = \{ x_1, \beta \}$, it follows $z_1 \in {\Cal U}_3
\cup {\Cal U}_4$. Let $z_1 \in {\Cal U}_3$. We choose $z_j \in \{ z_2, 
z_3 \} \setminus {\Cal U}_2$ and define $f(u_1^*) = \{ 
z_j \}$. By \C3, $u_1^* \in S_1 \cup S_6$.
Since $z_1 \in {\Cal U}_3 \cap N(z_j)$, $f(u_1^*)$ satisfies (7). 
If $z_1 \in {\Cal U}_4$, then let $f(u_1^*) = \emptyset$. 
Since $z_1 \in N(\beta)$,
by (8),  
$$|U_2^+| - |S_0| - |S_3| \geq - 2, ~\text{and if} ~(9) 
~\text{fails, then there is} ~z_1 \in {\Cal U}_4 ~\text{with} 
~|N(z_1) \cap (V \setminus {\Cal U}_4)| \geq 4. \tag11$$ 
\hskip 12pt
Next we compare $\sum \w(z)$ for $z \in U_3^2$ with
$|S_4|$. 
Denote $\phi(z) = |N(z) \cap {\Cal U}_2| +
0.5|N(z) \cap {\Cal U}_3|$. 
If $z \in U_3^2$, then $\phi(z) \geq 4$ and $\w(z) = \phi(z) - 2 \geq
\phi(z) - 0.5 \phi(z) \geq 
0.5|N(z) \cap {\Cal U}_2|$. Hence 
$$2\sum \{ \w(z): z \in U_3^2 \} \geq \sum \{ |N(z) \cap {\Cal U}_2|:
z \in U_3^2 \} \geq
|\{ y \in {\Cal U}_2: N(y) \cap U_3^2 \neq        
\emptyset \}| \geq |S_4|. \tag12$$
Recall $3 \leq |N(z) \cap {\Cal U}_2| + 0.5|N(z) \cap {\Cal U}_3| 
\leq 3.5$ for $z \in U_3^1$.
Since the vertices in $S_1$ satisfy (7), if $z \in U_3^1$ and $|\{ b 
\in S_1: z \in f(b) \}| = 3$, then $\w(z) = 1.5$.
We partition $U_3^1$ into six sets: \newline
Let $U_{3i}^1 = \{ z \in U_3^1: |\{ b \in S_1: z \in f(b) \}|
= i, ~\w(z) = 1.5 \}$, $i = 1, 2, 3,$ the set
$U_{3k}^1 = \{ z \in U_3^1: |\{ b \in S_1: z \in f(b) \}|
= k - 3, ~\w(z) = 1 \}$, $k = 4, 5$, and
$U_{36}^1 = U_3^1 \setminus (U_{31}^1 \cup ... \cup 
U_{35}^1).$
Thus $$|S_1| = |U_{31}^1| + |U_{34}^1| + 2(|U_{32}^1| +
|U_{35}^1|) + 3|U_{33}^1|. \tag13$$
Since each $y \in S_2$ is adjacent to two vertices in            
$f(y) \cap U_3^1$, and the vertices in $S_1$ satisfy (7),
it follows  
$2|S_2| \leq 2|U_{31}^1| + |U_{32}^1| + |U_{34}^1| + 
3|U_{36}^1|$, and by (13), 
$$\sum\{ \w(z): z \in U_3^1 \} \geq 1.5(|U_{31}^1| + |U_{32}^1| + 
|U_{33}^1|) + |U_{34}^1| + |U_{35}^1| + |U_{36}^1|  
\geq 0.5(|S_1| + |S_2|). \tag14$$
\hskip 12pt 
Let $w \in {\Cal U}_4$. Since $w \in N(\beta)$, the vertex $\beta \in 
N(x_k) \cap N(w)$, 
and $G$ is $K_{2,3}$-free, it follows
$|N(w) \cap {\Cal U}_2 \cap N(x_k)| \leq 1$, $k = 1, 2$, and
$|N(w) \cap S_5| \leq 2$. 
We partition ${\Cal U}_4$ into three sets. 
Let $U_{42} = \{ w \in {\Cal U}_4 \cap f(y) \cap f(b):
y, b \in S_5, y \neq b \}$, $U_{41} = \{ w \in ({\Cal U}_4 \cap 
f(y)) \setminus U_{42}: y \in S_5 \},$
and $U_4^* = {\Cal U}_4 \setminus (U_{41} \cup U_{42}).$ Thus 
$|S_5| \leq 2|U_{42}| + |U_{41}|$. 
If $b \in {\Cal U}_3$, then by \C3, $\w(b) \geq 0$. 
Since $f(\gamma)$ is defined when $u_1^*$ is the only vertex in 
${\Cal U}_2^-$ whose $f(u_1^*)$ is undefined, it follows
$|{\Cal U}_2^-| = |S_0| + ... + |S_6|.$
Note that $S_6 \subseteq \{ u_1^* \}$, $|S_6| \leq 1$, and 
$|S_6| \leq |U_3^0|$.
By (6), (12), and (14),
$$e(G) = e(V_1) + e([V_1, {\Cal U}_2]) + e({\Cal U}_2)
+ e([{\Cal U}_2, {\Cal U}_3]) + e({\Cal U}_3) + e([V \setminus 
{\Cal U}_4, {\Cal U}_4]) + e({\Cal U}_4)$$
$$\geq 2|V_1| - 3 + 2|{\Cal U}_2| + 0.5(|U_2^+| - \sum_{0 \leq i
\leq 6}|S_i|) + 2|{\Cal U}_3| + \sum_{z \in U_3^1 \cup U_3^0} \w(z)
+ 0.5|S_4| + 2|{\Cal U}_4| + \sum_{w \in {\Cal U}_4} \w(w)$$
$$\geq 2n - 3 + 0.5(|U_2^+| - |S_0| - |S_3|)
+ \sum_{w \in U_{42}} (\w(w) - 1) + \sum_{w \in U_{41}}
(\w(w) - 0.5) + \sum_{w \in U_4^*} \w(w). \tag15$$
If $w \in U_{42}$, then by \C3, $w$ is adjacent to  
three vertices in $V \setminus {\Cal U}_4$ and $\w(w) \geq 1$.
If $w \in U_{41}$, then 
$|N(w) \cap (V \setminus {\Cal U}_4)| \geq 2$,
and since $f(y)$ satisfies (7) for $y \in {\Cal U}_2^- \cup \{ \gamma 
\}$, it follows $d(w) \geq 3$ and $\w(w) \geq 0.5$.
Let $W = \{ w \in U_4^*: \w(w) < 0 \}$ and $w \in W$ with
$N(w) = \{ \beta, u \}$. 
We define $p(w) = u$, the neighbor of $w$ in ${\Cal U}_4$.
By our assumption, there is a unique vertex $z$ such
that $N(\beta) \cap N(u) = \{ w, z \}$.
We define $h(w) = z$. 
We define
$$W_1 = \{ w \in W: h(w) \in V \setminus {\Cal U}_4 \}, ~\text{the set}
~P = \{ p(w): w \in W_1 \},$$ $$W_2 = W \setminus W_1 = \{ w \in W: h(w) 
\in {\Cal U}_4 \}, ~\text{and} ~H = \{ h(w): w \in W_2 \}.$$
Let $u \in P$. Thus there is $w \in W_1$ such that $u = p(w)$. 
Since $N(u) \cap 
N(\beta) = \{ w, h(w) \}$ and $h(w) \in {\Cal U}_2$,
it follows $\w(u) \geq 0.5$, the vertex $h(w)$ is not adjacent to any 
vertex in 
$W$, $u \not \in W \cup H$, and $u \neq p(w')$ for any $w' \in W 
\setminus \{ w \}$. Thus $P \cap (W \cup H) = \emptyset$ and 
$|W_1| = |P|$. Let $u \in P \cap (U_{42} \cup U_{41})$. Thus there is
$y \in {\Cal U}_2^- \cup \{ \gamma \}$ such that $u \in f(y)$. If
$y = \gamma$, then by (10),
$|N(u) \cap (V \setminus {\Cal U}_4)| \geq 3$.
If $y \in {\Cal U}_2^-$, then $y \not \in N(\beta)$ and $y, h(w), \beta
\in N(u) \cap (V \setminus {\Cal U}_4)$. Since $w \in N(u) \cap {\Cal 
U}_4$, it follows $\w(u) - 1 \geq 0.5$. Let $z \in H$. Thus there is
$w \in W_2$ such that $z = h(w)$ and $p(w) \in N(z) \cap N(\beta)$.
Since $|N(z) \cap N(\beta)| \leq 2$, it follows $z = h(w)$ for at 
most two $w$ in $W_2$. Thus $|W_2| \leq 2|H|$.
Since $w \in A$, by our choice of $\alpha$, $d(z) \geq d(\beta) \geq 5$,
$\w(z) \geq 1$ and $z \not \in W$.
If $z \in U_{41}$, then $z$ is adjacent to two vertices in
$V \setminus {\Cal U}_4$, and $\w(z) - 0.5 \geq 1$.
If $z \in U_{42}$, then $z$ is adjacent to three vertices 
in $V \setminus {\Cal U}_4$, and $\w(z) - 1 \geq 1$.
If (9) holds, then by (15), $e(G) \geq 2n - 3.5 - 0.5|W_1| - 
0.5|W_2| + 0.5|P| + |H| \geq 2n - 3.5.$

Suppose (9) does not hold. By (11), there is $z_1 \in {\Cal U}_4$
with $|N(z_1) \cap (V \setminus {\Cal U}_4)| \geq 4$ and $\w(z_1)
- 1 \geq 1$. Recall $P \cap H = \emptyset$. If $z_1 \in H$, then 
since $d(z_1) \geq 5$, we have $\w(z_1) - 1 \geq 1.5$.
By (11) and (15), $e(G) \geq 2n - 4 - 0.5|W| + 0.5|P \setminus \{ z_1 \}|
+ |H \setminus \{ z_1 \}| + (\w(z_1) - 1) \geq 2n - 3.5$.
\qed\enddemo

\proclaim{Lemma 5.3} If $G$ is a sat$(n, K_{2,3})$-graph
with $\delta(G) = 2$, then $e(G) \geq 2n - 3$.
\endproclaim
\demo{Proof} By Lemma 5.2, we assume that degree $2$ vertices have 
nonadjacent neighbors.
Let $B$ be the set of degree $2$ vertices whose neighbors have
exactly one common neighbor.
Let $\alpha \in B$, if $B \neq \emptyset$. Let $N(\alpha) = \{
x_1, x_2 \}$. If $|B| = 0$, then $N(x_1) \cap N(x_2) = \{ \alpha, \beta 
\}$. If $\alpha \in B$, then we define $V_1 = N[\alpha]$, the set
$U_2 = (N(x_1) \cup 
N(x_2)) \setminus V_1$, and $U_3 = V \setminus (V_1 \cup U_2)$, otherwise,
let $V_1 = N[\alpha] \cup \{ \beta \}$, $U_2 = (N(x_1) \cup N(x_2)) 
\setminus V_1$, the set $U_3 = V \setminus (V_1 \cup U_2 \cup N(\beta))$,
and $U_4 = V \setminus (V_1 \cup U_2 \cup U_3).$ 
For $w \in U_4$, let $\w(w) = |N(w) \cap (V \setminus U_4)| + 0.5|N(w) 
\cap U_4| - 2$,
$$\ell = \sum_{w \in U_4} \w(w), 
~\theta_2 = \sum_{y \in U_2} (|N(y) \cap V_1| + 0.5|N(y) \cap 
U_2| - 2), ~\text{and} ~\theta_3 = \sum_{z \in U_3} (|N(z) \cap 
U_2| - 2).$$
Hence $$e(G) = e(V_1) + e([V_1, U_2]) + e(U_2)
+ e([U_2, U_3]) + e(U_3) + e([V \setminus U_4, U_4]) +
e(U_4)$$
$$= 2|V_1| - 4 + 2|U_2| + \theta_2 + 2|U_3| + \theta_3
+ e(U_3) + 2|U_4| + \ell
= 2n - 4 + \theta_2 + \theta_3 + e(U_3) + \ell.$$

We shall prove $\theta_2 + \theta_3 + e(U_3) + \ell \geq 0.5.$
By \C3, $\theta_2, \theta_3 \geq 0$.
First let $\alpha \in B$. Thus 
$N(x_1) \cap N(x_2) = \{ \alpha \}$ and $|U_4| = \ell = 0$.
Since $x_1 x_2 \not \in E(G)$, by \P1, without loss of            
generality, there is $z \in N(x_2)$ adjacent to 
$y_1, y_2 \in N(x_1)$. If $|N(y_i) \cap U_2| \geq 3$, then $\theta_2 
\geq 0.5$. By \C2, $|N(y) \cap N(x_i)| \geq 2$ for $y \in U_2$,
$i = 1$ or $2$. Thus let $N(y_i) \cap U_2 = \{ z, z_i \}$, where
$z_i \in N(x_2)$, $i = 1, 2$. Since $y_1 y_2 \not \in E(G)$, by \P1, 
there is $b \in N(y_k) \setminus \{ x_1 \}$ adjacent to 
$b_1, b_2 \in N(y_{3-k})$, $k = 1$ or $2$.
Since $d(b) \geq 3$, if $b \in U_3$, then $\theta_3 + e(U_3) \geq 1$.
Thus $b \in \{ z, z_k \}$. Since $y_{3-k} \in N(x_1) \cap N(b_i)$ and
$b \in N(x_2) \cap N(b_i)$, by \C3, if $b_i \in U_3$, then
$\theta_3 \geq 1$. Thus $\{ b_1, b_2 \} = \{ z, z_{3-k} \}$. 
Hence $b = z_k$, and $y_k, z, z_{3-k} \in N(z_k)$, and $\theta_2 \geq 0.5$.

Finally, let $B = \emptyset$. 
Since degree $2$ vertices have nonadjacent neighbors,
if $w \in U_4$ has $N(w) = \{ \beta,
u \}$, then $u \in V \setminus U_4$ and $\w(w) \geq 0$. 
Hence $\ell \geq 0$. Since
$\beta \not \in N(\alpha)$, by \C2, $x_i$ is adjacent to $b_1,
b_2$ in $N(\beta)$, $i = 1$ or $2$. 
If $b_1 \in N(b_2)$, then $\theta_2 \geq 1$.
Let $b_1 \not \in N(b_2)$. By \P1, there is $y \in N(b_j)$ 
adjacent to $y_1, y_2 \in N(b_{3-j})$, $j = 1$ or $2$.  
If $z \in \{ y, y_1, y_2 \}$ is in $U_4 \cap N(b_k)$, 
then since degree $2$ vertices have nonadjacent neighbors,
$N(z) \supsetneq \{ \beta, b_k \}$ and $\ell \geq 0.5$.
Let $\{ y, y_1, y_2 \} \cap U_4 = \emptyset$.
If one of $y_1, y_2, y$ is in $U_2 \cup \{ x_i \}$,
then $\theta_2 \geq 0.5$.
Hence $y, y_1, y_2 \in U_3$ and $e(U_3) \geq 2$.
\qed\enddemo  

\bigskip
\subhead \S\ 6. The case when $\delta(G) = 3$
\endsubhead
\medskip
In this section, we prove Theorem 1 when $\delta(G) = 3$. 
We define $\lambda(G) = \min \{ e(N(\alpha)): d(\alpha) = 3
 \}$ and discuss $e(G)$ in cases depending on $\lambda(G)$.

\proclaim{Lemma 6.1} If $G$ is a sat$(n, K_{2,3})$-graph
with $\delta(G) = 3$ and $\lambda(G) = 3$, then $e(G) \geq 2n -
2$.
\endproclaim
\demo{Proof} Let $R = \{ v \in V: d(v) = 3 \}$. 
We claim that there is $\beta \in V$ such that $R \subseteq N(\beta)$. 
Let $\alpha_1 \in R$ and $\alpha_2 \in R \setminus
N[\alpha_1]$. If $N(\alpha_1) = \{ b_1, b_2,
b_3 \}$, then since $G$ is $K_{2,3}$-free, it follows $N(b_1) \cap 
N(b_2) = \{ \alpha_1, b_3 \}$. Thus $|N(\alpha_1) \cap N(\alpha_2)| 
\leq 1$. We first show $|N(\alpha_1) \cap N(\alpha_2)| = 1$. By \P1, 
there exists $\beta \in N(\alpha_{3-i})$ adjacent to $y_1, y_2 \in
N(\alpha_i)$, $i = 1$ or $2$. We denote $N(\alpha_i)
= \{ y, y_1, y_2 \}$. Since $G$ is $K_{2,3}$-free and
$N(\alpha_i)$ is a clique, it follows $N(y_1) \cap
N(y_2) = \{ \alpha_i, y \}$. Thus $\beta = y$ and $\{ \beta \} =
N(\alpha_1) \cap N(\alpha_2)$. Let $\alpha_3 \in R \setminus N(\beta)$. 
Thus $\alpha_3 \not \in N(\alpha_i)$ and $N(\alpha_3)
\cap (N(\alpha_i) \setminus N(\alpha_{3-i})) = \{ z_i \}$, $i = 1,
2$. We denote $N(\alpha_3) = \{ z_1, z_2, z \}$. However $\beta,
\alpha_3, z \in N(z_1) \cap N(z_2)$, a
contradiction. Hence $R \subseteq N(\beta)$ and $d(\beta) \geq |R|$. 
Therefore, $$2e(G) = \sum \{ d(x): x \in V \}
\geq 3|R| + d(\beta) + 4(n - |R| - 1) \geq 4n - 4. \qed$$
\enddemo

\proclaim{Lemma 6.2} If $G$ is a sat$(n, K_{2,3})$-graph with 
$\delta(G) = 3$, then $e(G) \geq 2n - 3$.
\endproclaim
\demo{Proof} We choose a degree $3$ vertex $\alpha$ with
$e(N(\alpha)) = \lambda(G)$. By Lemma 6.1, we assume $\lambda(G)
\leq 2$.
Denote $N(\alpha) =
\{ x_1, x_2, x_3 \}$, where $x_3 \not \in N(x_2)$.
We define $V_0 = N[\alpha]$, ${\Cal U}_1 = \{ y \in V \setminus 
V_0: |N(y) \cap N(\alpha)| = 2 \}$, ${\Cal U}_2 = \{ y \in V 
\setminus V_0: |N(y) \cap N(\alpha)| = 1 \}$, ${\Cal U}_3 = V 
\setminus (V_0 \cup {\Cal U}_1 \cup {\Cal U}_2),$
$$\w(y) = |N(y) \cap (V_0 \cup {\Cal U}_1 \cup {\Cal 
U}_{i-1})| + 0.5|N(y) \cap {\Cal U}_i| - 2 ~\text{for} ~y \in 
{\Cal U}_i, ~\text{and} ~\theta_i = \sum_{y
\in {\Cal U}_i} \w(y), ~i = 2, 3.$$
Hence $$e(G) = e(V_0) + e([V_0, {\Cal U}_1]) + e({\Cal 
U}_1) + e([V_0 \cup {\Cal U}_1, {\Cal U}_2]) + e({\Cal U}_2) + 
e([V \setminus {\Cal U}_3, {\Cal U}_3]) + e({\Cal U}_3)$$
$$= 3 + \lambda(G) + 2|{\Cal U}_1| + e({\Cal U}_1) + 2|{\Cal 
U}_2| + \theta_2 + 2|{\Cal U}_3| + \theta_3
= 2n - 5 + \lambda(G) + e({\Cal U}_1) + \theta_2 + \theta_3.$$ 
It is sufficient to prove
$$\lambda(G) + e({\Cal U}_1) +
\theta_2 + \theta_3 \geq 1.5. \tag16$$

\noindent {\bf Claim 6.2.1.} If $\lambda(G) = 0$, then (16) holds.

\noindent {\it Proof.} Since $e(N(\alpha)) = 0$, by \C3, if $y
\in {\Cal U}_2 \cup {\Cal U}_3$, then $\w(y) \geq 0$. 
Let $W_k = \{ y \in {\Cal U}_2 \cup {\Cal U}_3: \w(y) \geq 0.5k
\}$. Note that $W_3 \subseteq W_2 \subseteq W_1$ and $\theta_2 + 
\theta_3 \geq 0.5(|W_1| + |W_2| + |W_3|)$. 
First let ${\Cal U}_1 = \emptyset$. Let there exist $b_i \in {\Cal 
U}_2 \cap N(x_i)$ adjacent to $z_i \in {\Cal U}_3$ or having 
$\w(b_i) > 0$ for $i = 1, 2$, and $3$. 
Since $\delta(G) = 3$, by \C3, $z_i \in W_1$, and if $z_i = z_j$ 
and $i \neq j$ then $z_i \in W_2$. Also, if $z_1 \in N(b_1) \cap
N(b_2) \cap N(b_3)$, then
$z_1 \in W_3$. Thus $|W_1| + |W_2| + |W_3| \geq 3$.
Hence we assume that if $y \in N(x_1)$, then $\w(y) = 0$ and 
$N(y) \subseteq V \setminus {\Cal U}_3$.
Since $\delta(G) = 3$, there are $y, y' \in N(x_1) \cap {\Cal U}_2$. 
By \C2, there are $y_1, y_2 \in N(x_i)$ such that $N(y) = \{ x_1, 
y_1, y_2 \}$. Without loss of generality, $i \in \{ 1, 2 \}$.
Since $N(x_1) \cap N(x_3) = \{ \alpha \}$ and $N(y_1) \cap N(y_2) = 
\{ x_i, y \}$, it follows $|N(b) \cap N(y)| \leq 1$ for all $b \in 
N(x_3)$. Since $y x_3 \not \in E(G)$, by \P1, $y_1$ or $y_2$, say $y_1$, 
is adjacent to $z_1, z_2 \in N(x_3)$. Thus $\w(y_1) \geq 0.5$. By our 
assumption, $i = 2$. Since $N(x_1) \cap N(x_2) = \{ \alpha \}$,
$N(y_1) \cap N(y_2) = \{ x_2, y \}$, and $y 
x_2 \not \in
E(G)$, by \P1, $y_j$ is adjacent to $v_1, v_2 \in N(x_2)$, $j = 1$
or $2$. If $j = 1$, then $y_1 \in W_3$. Thus $j = 2$ and $y_1, y_2
\in W_1$. If $y_k$ is adjacent to $y'$, $k = 1$ or $2$, then 
$y_k \in W_2$ and $|W_1| + |W_2| \geq 3$. Let $|\{ y_1, y_2 \} \cap
N(y')| = 0$. Similarly to $y$, the vertex $y'$ has a neighbor
$y_3 \in W_1$. Thus $|W_1| \geq 3$. Hence we assume ${\Cal 
U}_1 \neq \emptyset$. 

Let $u_1 \in N(x_1) \cap N(x_2)$ and let $|N(u_1) \cap {\Cal
U}_1| = 0$. By \C2, $u_1$ is adjacent to $y_1, y_2 \in 
N(x_i) \cap {\Cal U}_2$, where $i \in \{ 1, 3 \}$. Let $z_j 
\in N(y_j) \setminus \{ x_i, u_1 \}$, $j = 1, 2$. 
By \C3, $y_j \in W_2$, $y_j \in W_1$, or $z_j \in W_1$ when 
$z_j$ is in ${\Cal U}_1$, ${\Cal U}_2$, or ${\Cal U}_3$, 
respectively. Since $N(y_1) \cap N(y_2) = \{ x_i, u_1 \}$, 
it follows $z_j \not \in N(y_{3-j})$ and $|W_1 \cap \{ y_1, 
y_2, z_1, z_2 \}| \geq 2$.
If $e({\Cal U}_1) > 0$, or $z_j \in {\Cal U}_1$, or there is 
$z \in N(y_j) \setminus \{ u_1, x_i, z_j \}$, then (16) holds.
Thus $e({\Cal U}_1) = 0$ and $N(y_j) = \{ u_1, x_i, z_j \}$,
where $z_j \not \in {\Cal U}_1$, $j = 1, 2$.
Let there be $u_2 \in {\Cal U}_1 \setminus \{ u_1 \}$. 
By \C3, there are $y_3, y_4 \in N(x_k) \cap N(u_2)$. 
Similarly to $y_1, y_2$, $N(y_3) = \{ u_2, x_k, z_3 \}$ and 
$z_3 \not \in {\Cal U}_1$.
If $z_3 \in \{ y_1, y_2, z_1, z_2 \} \cap {\Cal U}_2$ or 
$z_3 \not \in \{ z_1, z_2 \}$, then $|W_1| \geq 3$. Thus 
$z_1 = z_3$ and $z_1 \in {\Cal U}_3$. Since $y_1 y_3 \not 
\in E(G)$ and $z_1 \in {\Cal U}_3$, by \P1, $z_1 \in
(N(u_1) \cup N(u_2)) \cap W_2$ and (16) holds. 
Hence ${\Cal U}_1 = \{ u_1 \}$. 

Since $d(x_2) \geq 3$, there is $y_5 \in N(x_2) \cap 
{\Cal U}_2$. 
Let $y_5 \in N(u_1)$ and $z_5 \in N(y_5) \setminus \{ u_1, 
x_2 \}$. If $y_5 = z_j$, $j = 1$ or $2$, then $y_j, z_j \in 
W_1$ and $|W_1| \geq 3$. Since $i \neq 2$, we assume $y_5 
\not \in \{ y_1, y_2, z_1, z_2 \}$. By \C3, $y_5 \in W_1$ or 
$z_5 \in W_1$ when $z_5$ is in ${\Cal U}_2$ or ${\Cal U}_3$,
respectively. Also, if $z_5 \in \{ z_1, z_2 \}$, then
$z_5 \in W_2$. Thus $|W_1| + |W_2| \geq 3$. Hence $y_5 \not 
\in N(u_1)$. By \C2, there are $y_6, y_7 \in N(x_k) \cap 
N(y_5)$. If there is $z \in N(y_5) \setminus \{ x_2, y_6, 
y_7 \}$, then by \C3, $\{ y_5, z \} \cap W_1 \neq \emptyset$.
Also, if $z \in \{ z_1, z_2 \}$, then $z \in W_2$.
Hence $N(y_5) = \{ x_2, y_6, y_7 \}$. 
Since $N(x_3) \cap N(x_2) = \{ \alpha \}$, $N(y_6) \cap 
N(y_7) = \{ x_k, y_5 \}$, $x_k \not \in N(x_3)$,
and $y_5 x_3 \not \in E(G)$, by \P1, $y_6$ or $y_7$, say 
$y_6$, is adjacent to $w_3, w_4 \in N(x_3)$. Thus $y_6 \in 
W_1$. Since $N(y_j) = \{ u_1, x_i, z_j \}$ and $d(y_6) \geq
4$, it follows $y_6 \neq y_j$, $j = 1, 2$. If 
$y_6 = z_j$, $j = 1$ or $2$, then $y_j \in W_1$. Thus $|W_1|
+ |W_2| \geq 3$. Hence $|N(u) \cap {\Cal U}_1| \geq 1$ if $u 
\in {\Cal U}_1$. 

We can assume $E(G[{\Cal U}_1]) = \{ u_1 u_2 \}$, $|{\Cal 
U}_1| = 2$, and $u_2 \in N(x_2) \cap N(x_3)$. Since $d(x_1) 
\geq 3$, there is $y^* \in N(x_1) \cap {\Cal U}_2$. If there
is $z \in ({\Cal U}_3 \cup {\Cal U}_1) \cap N(y^*)$, then 
since $\delta(G) = 3$, by \C3, $(\{ y^* \} \cup (N(y^*)
\cap {\Cal U}_3)) \cap W_1 \neq \emptyset$. 
Let $N(y^*) \subseteq \{ x_1 \} \cup {\Cal U}_2$.
If there is $v \in N(y^*)$ adjacent to distinct $v_1, v_2 
\in N(x_3)$, then $v \in {\Cal U}_2 \cap W_1$.
Let $|N(v) \cap N(x_3)| \leq 1$ for all $v \in N(y^*)$. 
Since $y^* x_3 \not \in E(G)$ and $N(x_3) \cap N(x_1) = \{ 
\alpha \}$, by \P1, there is $z \in N(x_3)$ adjacent 
to $z_1, z_2 \in N(y^*) \setminus \{ x_1 \}$. By \C2, $z_1 
\in W_1$ and (16) holds. \qed

From now on, we assume $\lambda(G) = 1$ or $2$ and $x_1 \in 
N(x_2)$.
Recall $x_3 \not \in N(x_2)$. We define
$$U_2^+ = \{ y \in {\Cal U}_2: \w(y) \geq 0.5 \},
~{\Cal U}_2^- = \{ y \in {\Cal U}_2: \w(y) < 0 \}, ~\text{and} ~{\Cal 
U}_{2i}^- = {\Cal U}_2^- \cap N(x_i), i = 1, 2, 3,$$
$${\Cal U}_3^g = \{ z \in {\Cal U}_3: d(z) \geq 4 ~\text{or}
~|N(z) \cap (V \setminus {\Cal U}_3)| \geq 3  
~\text{or} ~N(z) \cap ({\Cal U}_1
\cup ({\Cal U}_2 \setminus {\Cal U}_2^-)) \neq \emptyset \},$$
and ${\Cal U}_3^b = {\Cal U}_3 \setminus {\Cal U}_3^g$. For $y \in {\Cal 
U}_2^-$, we define $f(y)$, a subset of $N(y) \cap {\Cal U}_3$
satisfying (17):
$$\text{If} ~z \in f(y) ~\text{and} ~z \not \in {\Cal 
U}_3^g, ~\text{then there is} ~y_1 \in {\Cal U}_2^-
\cap N(z) ~\text{with} ~z \not \in f(y_1). \tag17$$
Let $y \in {\Cal U}_{2j}^-$ and $\w(y) = - 0.5$.
By \C3, there is $y' \in N(x_i)$ such that $N(y) \cap 
(V \setminus {\Cal U}_3) = N(y) \cap N(x_i) = \{ x_j, y' \}$. 
Since $d(y) \geq 3$, there is $z \in N(y) 
\cap {\Cal U}_3$. We choose one $z \in N(y) \cap {\Cal 
U}_3^g$ and define $f(y) = \{ z \}$, if $N(y) \cap {\Cal 
U}_3^g \neq \emptyset$. 
Let $N(y) \cap {\Cal U}_3^g = \emptyset$ and $z \in N(y) 
\cap {\Cal U}_3$. 
By \C2, there is $y_1 \in {\Cal U}_{2j}^-$ and $w \in {\Cal 
U}_3$ such that $N(z) = \{ y, y_1, w \}$. By \C3, $y y_1 
\not \in E(G)$. 
Since $N(y) \cap N(y_1) = \{ x_j, z \}$ and $\lambda(G) \geq 1$, it 
follows $w$ is adjacent to exactly one of 
$y, y_1$, $\lambda(G) = 1$, $\{ i, j \} = \{ 1, 2 \}$, and $x_3 \not \in 
N(x_j)$.
Since $N(y) \cap N(y_1) = \{ x_j, z \}$, $N(b) \cap N(z) \subseteq \{ 
w \}$ for 
all $b \in N(x_3)$, and $z x_3 \not \in E(G)$, by \P1, $w$ is adjacent to
$y_2, y_3 \in N(x_3)$. 
Thus $w \in {\Cal U}_3^g \cap N(y_1)$
and $z \not \in f(y_1)$. We define $f(y) = \{ z 
\}$ and $f(y)$ satisfies (17).

Next let $y \in {\Cal U}_2^-$ with $\w(y) = - 1$. By \C3, $y \in N(x_1)$,
$x_1 \in N(x_2) \cap N(x_3)$, and $\lambda(G) = 2$.
Since $d(y) \geq 3$ and $N(y) \cap (V \setminus {\Cal U}_3) = \{ x_1 \}$,
there are $z, z' \in N(y) \cap {\Cal U}_3$. We claim $z, z' \in
{\Cal U}_3^g$. Let $z \in {\Cal U}_3^b$ and $N(z) = \{  y, y_1, w \}$,
where $y_1 \in {\Cal U}_2^-$ and $w \in {\Cal U}_3$.
By \C2, 
$y_1 \in N(x_1)$ and $y y_1 \not \in E(G)$. Since $\lambda(G) = 2$, it 
follows $x_1, z, w \in N(y) \cap N(y_1)$, a contradiction  
justifying our claim. 
We define $f(y) = \{ z, z' \}$ and $f(y) \subseteq
{\Cal U}_3^g$ satisfying (17). 
 
We claim that
$$\text{if} ~z \in {\Cal U}_3, ~\text{then} 
~\w(z) \geq \epsilon + 0.5|\{ y \in {\Cal U}_2^-: z \in f(y) \}|,
\tag18$$
where $\epsilon = 0.5$ if $x_3 \not \in N(x_1)$ and $|N(z) \cap N(x_3) 
\cap {\Cal U}_2| \geq 1$ or $|N(z) \cap N(x_3) \cap {\Cal U}_1| \geq 2$, 
and $\epsilon = 0$, otherwise.

Let $z \in {\Cal U}_3$. 
By \C2, $N(z) \cap ({\Cal U}_1 \cup {\Cal U}_2) \neq \emptyset$.
Also, if $x_3 \not \in N(x_1)$, then
$N(x_3) \cap {\Cal U}_2 \subseteq U_2^+$, since $x_3 \not \in N(x_2)$.
Since $d(z) \geq 3$, (18) follows \C3 if $z \not \in f(y)$ for any 
$y$. Let $z \in f(y)$, where $y \in {\Cal U}_{2j}^-$. 
By \C3, $|N(z) \cap (V \setminus {\Cal U}_3)| \geq 2$.
First let $N(z) \cap (V \setminus {\Cal U}_3) = \{ y, y_1 \}$.
Since $d(z) \geq 3$ and $f(y)$ satisfies (17), it follows $z \in {\Cal 
U}_3^g$ or $z \not \in f(y_1)$, and (18) holds.
Let $|N(z) \cap (V \setminus {\Cal U}_3)| = 3$. 
Since $\w(z) \geq 1 + 0.5(d(z) - 3)$, (18) 
holds if $d(z) \geq 4$ or $|N(z) \cap {\Cal U}_1| > 0$. 
Let $d(z) = 3$ and $N(z) \subseteq {\Cal U}_2$.
If $\lambda(G) = 2$, then $\epsilon = 0$, there is $y' \in
N(z) \setminus {\Cal U}_2^-$ with $|N(y') \cap N(z)| = 2$, and (18) holds.
Let $E(G[N(\alpha)]) = \{ x_1 x_2 \}$ and $y_1 y_2 \in
E(G[N(z)])$. If $y_i \not \in {\Cal U}_2^-$, $i = 1$ or $2$, then (18) 
holds. By \C3, $y_i \in N(x_i) \cap 
{\Cal U}_2^-$, $i = 1, 2$. Let $y_3 \in N(z) \setminus \{ y_1, y_2 \}$. 
Since $z \not \in N(\alpha)$, by \C2, $y_3 \in N(x_1) \cup N(x_2)$. 
Since $z x_3
\not \in E(G)$, by \P1, there is $b \in N(x_3) \cap N(y_3)$. 
By \C3, $y_3 \not \in {\Cal U}_2^-$ and (18) holds.
If $|N(z) \cap (V \setminus {\Cal U}_3)| \geq 4$,
then $\w(z) \geq 0.5|N(z) \cap (V \setminus {\Cal U}_3)|$ and (18) 
holds. 

By our construction, if $\w(y) = - 0.5$, then $|f(y)| = 1$. Also, if
$\w(y) = - 1$, then $|f(y)| = 2$. Let $\ell^* = \sum \{ \w(y): y \in 
U_2^+ \}$. By (18), 
$$\theta_2 + \theta_3 = \sum \{ \w(y): y \in U_2^+ \cup {\Cal U}_3 \cup 
{\Cal U}_2^- \}
\geq \ell^* + \epsilon + 
0.5\sum_{y \in {\Cal U}_2^-}|f(y)| 
+ \sum_{y \in {\Cal U}_2^-}\w(y) \geq \ell^* + \epsilon, \tag19$$
where $\epsilon = 0.5$ if $x_3 \not \in N(x_1)$ and there is $z \in {\Cal
U}_3$ with $|N(z) \cap N(x_3) \cap {\Cal U}_2| \geq 1$ or $|N(z) \cap 
N(x_3) \cap {\Cal U}_1| \geq 2$, and $\epsilon = 0$, otherwise.          

By (19), if $\lambda(G) = 2$ or $e({\Cal U}_1) \geq 1$, then
(16) holds. Let $\lambda(G) = 1$, $E(G[N(\alpha)]) = 
\{ x_1 x_2 \}$, and $e({\Cal U}_1) = 0$.
If there is $y \in N(x_3) \cap {\Cal U}_2$ adjacent to a 
vertex in ${\Cal U}_3$ or with $\w(y) \geq 0.5$, then by 
(19), (16) holds.
Thus we assume $\w(y) = |N(y) \cap {\Cal U}_3| = 0$ for
$y \in N(x_3) \cap {\Cal U}_2$.
Since $\delta(G) = 3$, $e({\Cal U}_1) = 0$, and $\lambda(G) 
= 1$, there is $y^* \in N(x_3) \cap {\Cal U}_2$. By \C2, 
there are $x_j \in N(\alpha)$ and $y_1, y_2 \in N(x_j)    
\setminus {\Cal U}_1$ such that $N(y^*) = \{ x_3, y_1, y_2 
\}$. Let $j \neq 3$. Since $\lambda(G) = 1$, $y_1 y_2 \in 
E(G)$. By \C2, either $N(y_1) \cap {\Cal U}_1 \neq 
\emptyset$ or there is $x_i \in N(\alpha)$ with $|N(y_1) 
\cap N(x_i)| \geq 2$. Thus $\ell^* \geq \w(y_1) \geq 0.5$ 
and (16) holds. 
Hence $y_1, y_2 \in N(x_3)$ and similarly for $y_i$,  
$N(y_i) \subseteq \{ x_3 \} \cup (N(x_3)
\cap {\Cal U}_2)$, $i = 1, 2$. 
First let $|{\Cal U}_3| = 0$. Since $N(b) \cap N(y^*) \subseteq \{ x_3 
\}$ for all $b \in N(x_1)$ and $y^* x_1 \not \in E(G)$, by \P1, there is
$z \in N(y^*)$ with $|N(z) \cap N(x_1)| \geq 2$. Since $N(y_i) \subseteq 
\{ x_3 \} \cup (N(x_3) \cap {\Cal U}_2)$, $i = 1, 2$, it follows $z = x_3$.
There is $u \in (N(x_3) \cap N(x_1)) \setminus \{ \alpha \}$ and $u \in 
{\Cal U}_1$.  
Since $\delta(G) = 3$ and $|{\Cal U}_3| = e({\Cal U}_1) = 0$, there is 
$y \in N(u) \cap {\Cal U}_2$, $\ell^* \geq \w(y) \geq 0.5$, and (16) 
holds. Hence 
there is $z^* \in {\Cal U}_3$. 
Since $N(b) \cap N(y^*) \subseteq \{ x_3 \}$ for all $b \in N(z^*)$ and 
$y^* z^* \not \in E(G)$, by \P1, there is $w \in N(y^*)$ with $|N(w) \cap
N(z^*)| \geq 2$. Since $N(y_i) \subseteq \{ x_3 \} \cup (N(x_3) \cap 
{\Cal U}_2)$, $i = 1, 2$, and $|N(y) \cap {\Cal U}_3| = 0$ for all $y \in
N(x_3) \cap {\Cal U}_2$, it follows $w = x_3$, $|N(x_3) \cap {\Cal U}_1 
\cap N(z^*)| \geq 2$, and $\epsilon = 0.5$.
By (19), (16) holds justifying Lemma 6.2. 
\qed\enddemo

\noindent {\bf Acknowledgments} The author is grateful to 
Zolt\'an F\"uredi, Ko-Wei Lih, 
and John McDonald for their helpful suggestions on the 
preparation of this paper.
She also thanks the very kind referees for their careful 
reading and critical, helpful suggestions.

\enddocument

\end

\ref \no
\by 
\paper 
\jour Discrete Math.
\vol 67
\yr 19
\pages --
\endref

\endRefs

\Refs
\parskip = 2pt plus 2pt
\widestnumber\key{ERS}

\ref\no 1
\by Y. Ashkenazi
\paper $C_3$ saturated graphs
\jour Discrete Math.
\vol 297
\yr 2005
\pages 152--158
\endref             

\ref\no 2              
\by C. A. Barefoot; L. H. Clark; R. C. Entringer;
T. D. Porter; L. A. Sz\'ekely and Zs. Tuza
\paper Cycle-saturated graphs of minimum size.
Selected papers in honour of Paul Erd\H{o}s on the
occasion of his 80th birthday (Keszthely, 1993)
\jour Discrete Math.
\vol 150
\yr 1996
\pages 31--48  
\endref             

\ref\no 3
\by T. Bohman, M.Fonoberova, and O. Pikhurko 
\paper The saturation function of complete partite graphs
\jour Journal of Combinatorics
\vol 1
\yr 2010
\pages 149--170
\endref

\ref\no 4
\by B. Bollob\'as
\paper On generalized graphs
\jour Acta Math. Acad. Sci. Hungar
\vol 16
\yr 1965
\pages 447--452                           
\endref

\ref\no 5
\by G. Chen; R. Faudree; and R. Gould
\paper Saturation numbers of books
\jour Electron. J. Combin. 
\vol 15
\yr 2008
\pages Research Paper 118, 12 pp
\endref

\ref\no 6                                 
\by Y. Chen
\paper Minimum $C_5$-saturated graphs     
\jour  J. Graph Theory                  
\vol 61             
\yr 2009                 
\pages 111--126                      
\endref

\ref\no 7       
\by Y. Chen      
\paper All Minimum $C_5$-saturated graphs
\jour  J. Graph Theory  
\vol 67
\yr 2011 
\pages 9--26
\endref

\ref \no 8
\by P. Erd\H{o}s, Z. F\"uredi and Zs. Tuza
\paper Saturated $r$-uniform hypergraphs
\jour Discrete Math.
\vol 98                  
\yr 1991
\pages 95--104   
\endref

\ref\no 9
\by P. Erd\H{o}s, A. Hajnal and J. W. Moon
\paper A problem in graph theory
\jour Amer. Math. Monthly
\vol 71 
\yr 1964      
\pages 1107--1110     
\endref

\ref\no 10
\by J. Faudree, R. Faudree, and J. Schmitt
\paper A Survey of Minimum Saturated Graphs
\jour Electron. J. Comb.
\yr 2011
\pages DS19, Dynamic Survey, 36 pp
\endref

\ref\no 11
\by D. C. Fisher; K. Fraughnaugh and L. Langley
\paper On $C_5$-saturated graphs with minimum size.
Proceedings of the Twenty-sixth
Southeastern International Conference on  
Combinatorics, Graph Theory and Computing (Boca
Raton, FL, 1995)                     
\jour Congr. Numer.
\vol 112
\yr 1995              
\pages 45--48
\endref

\ref \no 12      
\by Z. F\"uredi and Y. Kim                         
\paper Cycle-saturated graphs with minimum number of edges
\pages arXiv:1103.0067
\endref


\ref \no 13
\by R. Gould; T. Luczak and J. Schmitt
\paper Constructive upper bounds for cycle-saturated graphs
of minimum size
\jour Electron. J. Combin.
\vol 13
\yr 2006
\pages Research Paper 29, 19 pp
\endref      

\ref \no 14
\by R. Gould and J. Schmitt
\paper Minimum degree and the minimum size of $K_2^t$-saturated
graphs
\jour Discrete Math.                 
\vol 307
\yr 2007
\pages 1108--1114
\endref

\ref \no 15
\by L. K\'aszonyi and Zs. Tuza
\paper Saturated graphs with minimal number of edges
\jour  J. Graph Theory                
\vol 10
\yr 1986          
\pages 203--210           
\endref

\ref \no 16
\by L. T. Ollmann    
\paper $K_{2,2}$ saturated graphs with a minimal
number of edges
\inbook Proceedings of the Third Southeastern 
Conference on Combinatorics, Graph Theory, and      
Computing (Florida Atlantic Univ., Boca Raton,
Fla., 1972)       
\yr 1972
\pages 367--392
\publ Florida Atlantic Univ.
\publaddr Boca Raton, Fla.
\endref    

\ref \no 17
\by O. Pikhurko
\paper The minimum size of saturated hypergraphs
\jour Combin. Probab. Comput.                 
\vol 8
\yr 1999   
\pages 483--492
\endref        

\ref \no 18
\by O. Pikhurko
\paper Results and open problems on minimum     
saturated hypergraphs
\jour Ars Combin.            
\vol 72
\yr 2004
\pages 111--127
\endref               

\ref \no 19
\by O. Pikhurko and J. Schmitt
\paper A note on minimum $K_{2,3}$-saturated graphs
\jour Australas. J. Combin.
\vol 40
\yr 2008
\pages 211--215
\endref

\ref \no 20
\by M. Truszczy\'nski and Zs. Tuza
\paper Asymptotic results on saturated graphs
\jour Discrete Math. 
\vol 87
\yr 1991
\pages 309--314
\endref               

\ref \no 21
\by Zs. Tuza
\paper $C\sb 4$-saturated graphs of minimum size
\jour Acta Univ. Carolin. Math. Phys.
\vol 30
\yr 1989
\pages 161--167  
\endref      

\endRefs
\enddocument

\end